\input amstex \input amssym
\documentstyle{amsppt}
\magnification1200
\pagewidth{6.5 true in}
\pageheight{9  true in}
\hoffset=-0.1true in
\NoBlackBoxes

\topmatter
\title Gaps
between fractional parts, and additive combinatorics
\endtitle
\author Antal Balog,
Andrew Granville and Jozsef Solymosi
\endauthor
\abstract We give bounds on the number of distinct differences $N_a-a$ as $a$ varies over all elements of a given finite set $A$, and $N_a$ is a nearest neighbour to $a$.
\endabstract
\thanks{ The first author is partially support by Hungarian National Research Grant K 109789, the second author is partially supported by NSERC, and the third author is partially supported by ERC Advanced Research Grant
no 267165 (DISCONV), by Hungarian National Research Grant NK 104183, and by NSERC.}
\endthanks
\endtopmatter
\document

\head 1. Introduction \endhead

Let $\alpha$ be a real, irrational number and
consider the set of points
$$
S_\alpha(N)=\{ \alpha n:\ 1\leq n \leq N \}
\subseteq \Bbb R/\Bbb Z,
$$
which is isomorphic to the set of points $\{  e^{2i\pi \alpha n}:\ 1\leq n \leq N \}$
on the unit circle $S^1$.  Moreover, one can view $S_\alpha(N)$ as the set of points
$\{\{\alpha n\} : n \leq N \}\subseteq [0,1)$,
where $\{t\} $ denotes the fractional part of real number $t$ (that
is, $\{t\} = t-[t] $ where $[t]$ is the largest integer $ \leq t $), and then
order them as
$$S_\alpha(N) =\{0\leq b_1<b_2<\ldots<b_N<1\}.$$

In 1957 Steinhaus observed, and in 1958 Vera S\'os proved,
[7,8]\footnote{As well as Sur\'anyi, Swierczkowski [9],
Sz\"usz, and Erd\H os and
Tur\'an}, that there are, at most,
three distinct
consecutive differences in $S_\alpha(N)$, that is
$$
|\{b_{i+1} -b_i : 1 \leq i \leq N\}| \leq 3
$$
where we take
$ b_{N+1} = 1 + b_1 $.
This is arguably surprising since it is well known that
$S_\alpha (N)$ becomes increasingly uniformly distributed $ \
\text{mod}\ 1$ as $N$ gets larger, that is $S_\alpha(N) $ looks random
globally for large $N$, whereas S\'os's result tells us that $S_\alpha(N) $
looks highly structured locally. In 2002, V\^{a}j\^{a}itu and Zaharescu [10]
considered a question in-between these two extremes:
Take the set of $b_i$'s above and erase as many elements of the set as
one likes -- how large can one make the resulting set of differences?
They proved the following result (though with the constant
$2+2\sqrt{2}$ in place of our $2\sqrt{2}$):

\proclaim{Corollary 2} For any subset A 
of $S_\alpha(N) $ there are no more
than $2 \sqrt{2N} + 1 $ distinct consecutive
differences in A; that is, if
$A = \{ 0 \leq a_1 < a_2 <...<a_m <1\leq a_{m+1}:=1+a_1\}
\subset S_\alpha (N)$
then
$$
|\{a_{j+1} -a_j : 1\leq j \leq m \}| \leq 2
\sqrt{2N} + 1 .\tag 1
$$
Moreover there exists  ${A}\subseteq S_\alpha (N) $ with at least
$\sqrt{2N} -1$ distinct consecutive differences, so the result here
is ``best possible'' up to the value of
the constant in front of the ``$\sqrt{N}$''.
\endproclaim

The proof in [10] comes from careful number theoretic considerations, whereas
our proof stems  from more general additive combinatorics ideas, emerging from
the fact that $S_\alpha (N)+S_\alpha (N)$ is not much larger than $S_\alpha (N)$. Here and throughout we define
$$
A+B:=\{ n:\ n=a+b \  \text{for some}\ a\in A, b\in B\} .
$$

A set of points in $\Bbb R/\Bbb Z$ has a natural ``circular'' ordering by starting at any given point and then proceeding anticlockwise. Evidently there is no generic place to start the ordering but this will not effect on our results.  Each of the points will have a {\sl smaller} and a {\sl larger} neighbour, the two points that appear before and after our given point in the ordering, respectively (where, for the purposes of this definition, the first point in our ordering appears after the last).

We need a notion of length in $\Bbb R/\Bbb Z$: Any $t\in \Bbb R/\Bbb Z$ represents the coset $t+\Bbb Z$, and we define $\| t\|$ to be the minimal size  of any element of $t+\Bbb Z$, that is
 $$
\| t \|:= \min_{n\in \Bbb Z} |t-n|.
$$
Hence the distance between two points in $x,y\in  \Bbb R/\Bbb Z$ is given by $\| x-y\|$. 


\proclaim{Theorem 1} If $B$ is a finite subset of $\Bbb R/\Bbb Z$
then any subset $A$ of $B$ has at
most
$$
\sqrt{2|B|}\ \frac{|A+B|}{|B|} + 1
$$
distinct, consecutive differences. Moreover there are
subsets $A$ of $B$
with at least $\sqrt{2 |B|} -1 $ distinct consecutive differences.
\endproclaim

The  bound in Theorem 1 is better than the trivial bound, $\leq |A|$,
if $|A+B|\ll |A|\sqrt{|B|}$.

A map $\phi: A\to B$ is called a {\sl $2$-isomorphism} if $a_1-a_2=a_3-a_4$ if and only if $\phi(a_1)-\phi(a_2)=\phi(a_3)-\phi(a_4)$ for all $a_1,\ a_2,\ a_3,\ a_4$. Note that this does not effect coincidences amongst the set of differences of elements of $A$. Any non-degenerate affine transformation is a $2$-isomorphism on the reals, so that any given set of reals is $2$-isomorphic to a subset of $[0,\frac 12)$. The natural embedding $[0,\frac 12)\to \Bbb R/\Bbb Z$ given by sending $t$ to $t+\Bbb Z$ is  a $2$-isomorphism, and  the composition of these two $2$-isomorphisms is also a $2$-isomorphism. Hence we can apply Theorem 1 to any  finite  set of reals.

For any ordered set $b_1<b_2<\ldots <b_m$ define
$$D(B):=\{ b_{i+1}-b_i:\ 1\leq i\leq m-1\}$$
to be the set of distinct consecutive differences between elements of $B$.
Since $A+B\subseteq B+B$,
we immediately deduce the following result from Theorem 1:

\proclaim{Corollary 1} If $B$ is a finite set of real numbers then
$$
\sqrt{2 |B|} -1 \leq \max_{A\subseteq B} |D(A)| \leq \sqrt{2 |B|}
\cdot \frac{|B+B|}{|B|} + 1.
$$
\endproclaim

A nice example is given by the projection of the lattice points of a compact
subset of a $k$-dimensional lattice onto $\Bbb R / \Bbb Z$. That is
for given real numbers $\alpha_1,...,\alpha_n $ and integers $ N_1,...,N_k\geq 1$,
the set
$$
B =\left\{ n_1 \alpha _ 1 + n_2 \alpha _2 + ...+ n_k \alpha_k \ :\ {0 \leq n_j < N_j} \
\text{for}\ {1\leq j \leq k} \right\}    \subset \Bbb R/\Bbb Z \tag{2}
$$
typically satisfies $| B + B | \leq 2^k | B|$,
so the upper and lower bounds in  Corollary 1 are equal  up to a factor $\lesssim 2^k$.
Cobeli et al [3], discovered and proved the following, wonderful
generalization of Steinhaus's problem: \ If the elements of  $B$ are
$ 0 \leq b_1 < b_2 <...< b_N < 1 $
then there exists a set $\{ r_1,...,r_l\}$ of positive real numbers with $l\leq 2^k$,
such that every difference
$b_j - b_i $ with $ 1\leq i < j \leq N $ can
be written as a sum of the $r_i$'s. Rather more generally:

\proclaim{Theorem 2} Let $B$ be a finite subset of $\Bbb R/\Bbb Z$,
and $C$ be any subset of $B$ for which $ C-B = B-B $.
Let $\,b_c^-$ and $b_c^+$ be the smaller and larger neighbours of $c$ in
$B$, and let $R^-:=\{ c-b_c^-:\ c\in C\}$, $R^+:=\{ c-b_c^+:\ c\in C\}$.
Then every element of $B-B$ is a sum of elements of $R^-$, as well as a sum of elements of $R^+$,
that is $B-B \subset \langle R^- \rangle_{\Bbb N_0} = \langle R^+ \rangle_{\Bbb N_0}$. 
(Here $\Bbb N_0=\Bbb N\cup \{0\}$.)
\endproclaim

For $B$ as in (2) we can take
$$
C = \{ \delta_1 N_1 \alpha_1 + \delta_2 N_2 \alpha _2 +...+ \delta_k N_k
\alpha_k : \ \text{each} \
\delta_i =0 \ \text{or} \ 1 \}
$$
so that $ |R^-|,|R^+| \leq |C| \leq 2^k $, which is comparable to
$|B+B|/|B|$.

However, we have been unable, in general, to estimate the size of
$C$ in terms of  quantities associated with $B$, and there
does not seem  to be a close link between $|C|$ and $|B+B|/|B|$
in general. Let $r_3(N)$ denote the size of the largest set $S\subset \{1,\dots,N\}$ containing no non--trivial 3--term arithmetic progressions. As is well-known, $ r_3(N) \gg Ne^{-\delta\sqrt{\log N}}$ for some  constant $\delta>0$, thanks to a clever construction of Behrend [1].

\proclaim{Proposition 1} For any given integer  $N$ there exists a finite set of $N$
integers $B$ for which $|B+B|\leq 10N$ and if $C\subset B,\ C-B=B-B$ then $|C| \geq r_3(N)$.
\endproclaim

The length of $t\in (\Bbb R/\Bbb Z)^d$ is defined to be the minimal length  of any element of $t+\Bbb Z^d$; that is
 $$
\| t \|:= \min_{n\in \Bbb Z^d} |t-n| = \left(\sum_{i=1}^d \min_{n\in\Bbb Z}|t_i-n_i|^2\right)^{1/2}
$$
where $t=(t_1,\ldots,t_d)$. The distance between $x,y\in (\Bbb R/\Bbb Z)^d$ is then given by $\| x-y\|$.

For a given finite subset $A$ of $(\Bbb R/\Bbb Z)^d$ and $a\in A$,
let $N_a$ be one of the  elements of $A$ that is nearest to $a$ in the $\| .\|$-norm.. We now investigate the number of possible vectors $N_a-a$ as we vary over $a\in A$.

\proclaim{Theorem 3} Let $w_n=(\alpha_1n,\alpha_2n,\ldots,\alpha_dn) \in (\Bbb R/\Bbb Z)^d$.
For each $i\in [1,N]$ select $j=j(i)\in [1,N]$ with $j\ne i$ so that
$\|w_i-w_j\|$ is minimized. There  are   $\ll (4/3)^d$ distinct
elements of the set $D=\{ \pm (w_i-w_{j(i)}): i=1,2\ldots, N\}$.
\endproclaim

Unfortunately for a general set $A\in (\Bbb R/\Bbb Z)^d$ we can say somewhat less, namely we can prove a similar bound only for a larger subset of $A$.

\proclaim{Theorem 4} Let $A$ be a  finite subset of $(\Bbb R/\Bbb Z)^d$, and $\epsilon>0$ be a small real number. There exists $A'\subset A$, with $\geq (1-\epsilon)|A|$
elements, such that
$$
|\{ N_a-a:\ a\in A'\}|\ll (4/3)^d  \cdot \min_B \frac{|A+B|^2}{|A||B|} \ \cdot \  \frac 1{\epsilon} \log \frac 1{\epsilon}.
$$
\endproclaim

We give an example in section 5 to show that this cannot be much improved.

It would be interesting to determine $\min_B \frac{|A+B|^2}{|A||B|}$ for an arbitrary
set $A$ of real numbers. Taking $B=\{b\}$, any one element set, or taking $B=A$, evidently shows
$$\min_B \frac{|A+B|^2}{|A||B|} \leq \min\{ |A|, (|A+A|/|A|)^2\}.$$

\bigskip\head 2. The Three Gaps theorem and beyond \endhead

We begin with recalling Liang's elegant solution [6] to  Steinhaus's problem. We do so not only because it inspires some of our arguments, but also because we want to demonstrate how simple and short it is. 

\proclaim{The Three Gaps theorem}  Reorder $S_\alpha(N)$
as $0 \leq \{a_1\alpha\}<\{a_2\alpha\}<...<\{a_N\alpha\}< 1$.
Then the consecutive differences  $\{ a_{i+1}\alpha\} - \{a_i\alpha\},
1\leq i\leq N-1$ each equal  one of three numbers, namely the distances
between each pair of points amongst $\{a_N\alpha\}-1, 0 $ and $\{ a_1 \alpha\}$.
\endproclaim

\demo{Proof sketch}
$D(A)$ is contained in the set $G$ of gaps
$\{a_{i+1}\alpha\} - \{a_i\alpha\}$ for which $\{(a_i-1)\alpha\}$ and
$\{(a_{i+1}-1)\alpha\}$ are {\sl not} consecutive elements of the re-organized sequence
(that is, we select one representative from each set of gaps that
are obviously the same length.)  Note that if $\{(a_i-1)\alpha\}$ and
$\{(a_{i+1}-1)\alpha\}$ are {\sl not} consecutive elements of the re-organized sequence
then either $a_i=1$ or $a_{i+1}=1$, or there is some $\{a_j\alpha\}\in (\{(a_i-1)\alpha\}, \{(a_{i+1}-1)\alpha\})$. This last possibility implies that
$\{(a_j+1)\alpha\}\in (\{a_i\alpha\}, \{a_{i+1}\alpha\})$, and so $a_j=N$. Hence $|G|\leq 3$ and so  there are at most three distinct gaps. When $a_i=1$ the gap is the distance between
$0$ and $\{(a_{i+1}-1)\alpha\})$; when $a_{i+1}=1$  the gap is  the distance between
$\{(a_i-1)\alpha\}$ and $0$; and when
$\{a_i\alpha\}<\{ (N+1)\alpha\}< \{a_{i+1}\alpha\})$ the gap is the sum of the  two
previous distances.
\enddemo

\remark{Example} The Weierstrass parametrization
$\Bbb C/\Lambda \rightarrow E$ of an elliptic curve $E$, is an
isomorphism sending $ z\in \Bbb C $ to $P_z :=(\wp(z), \wp' (z))$, for a
certain function $\wp$ and $2$-dimensional lattice $\Lambda$. Hence
$P_{mz}=mP_z$. The
multiples $P, 2P,\ldots, NP$ of a given point $P$ on the real locus of $E$
(for example, the point $P =(1,3)$ on $E: y^2 = x^3 + 8x$) all lie
on the real locus of $E$, which is one continuous curve. Moving along
this curve from $-\infty$ to $+\infty$, we encounter these points in
some order $a_1 P, a_2 P,..., a_NP$. As noted in [5], the Three Gaps
theorem implies that there are, at most, three distinct differences
$(a_{i+ 1} - a_i)P$.
\endremark

Chung and Graham [2] came up with a beautiful generalization of the Three Gaps theorem, and Liang's proof works here too.

\proclaim{The $3k$-Gaps theorem}  Let $A$ be the set of elements $\beta_i+n_i\alpha \in \Bbb R/\Bbb Z $ for
$1\leq n_i\leq N_i$ and $i=1,2,\ldots,k$. Then the sizes of consecutive differences of $A$ have $\leq 3k$ different values.
\endproclaim

\demo{Proof}  List the elements of $A$ in an anticlockwise  ``circular order'', $a_1,a_2,\dots,a_N$, where each $a_j$ is in the form $\beta+n\alpha$ with some $\beta=\beta_i$ and $n=n_i$. Note that $N=N_1+\dots+N_k$, and $D(A)$, the set of consecutive differences is $\{a_2-a_1, a_3-a_2,\dots,a_N-a_{N-1},a_1-a_N\}$. We show that $|D(A)|\leq 3k$. 
$D(A)$ is contained in the set of gaps (better to call them arcs) $a_{i+1}-a_i$ for which $a_i+\alpha$ and $a_{i+1}+\alpha$ are not consecutive elements of $A$. This can only happen if either one of these two points (or both) is not an element of $A$, or there is an element $a\in A$ between them, but $a-\alpha$ is not an element of $A$. That is either one of the endpoints of the $A$--free arc $(a_i,a_{i+1})$ is of the form $\beta_j+N_j\alpha$ or there is a point of the form $\beta_j$ on the $A$--free arc $(a_i,a_{i+1})$. That is no more then 3 options for each $j$.
\enddemo

\bigskip\head 3. Bounding the number of differences in $A$
\endhead

\demo{Proof of the first part of Theorem 1}
Order the elements of $A$ as $a_1, a_2,\ldots, a_m$, going anti-clockwise around the circle. 
For each $d\in D(A)$ select $i(d) $
for which $ a_{i(d)+1} - a_{i(d)} = d $ and let
$$
J_A : = \{ i (d) : d \in D(A) \}
$$
For a $k$ chosen optimally later, we partition the set $ A+B $ into
${k}$ arcs $ I_1,I_2,...,I_k $ each
containing roughly the same number of elements of $
A+B $. We count the
number, $P $, of pairs $(i,b) $ with $ i \in J_A $ and $ b\in B $ for
which
${a_i + b} $ and $ {a_{i + 1} + b } $ both lie in the same arc
$I_j $.
For each fixed $b$, the points of $A+b$ follow each other in the order
$$
a_1 + b , a_2 + b , \dots , a_m + b
$$
and so consecutive numbers here lie in the same arc unless they
straddle the
boundary between consecutive arcs, which can happen for at
most $ k $ pairs.
Therefore
$$
P \geq |B| ( |D(A)| - k) \tag {3}
$$
On the other hand, if we are
given a pair of integers $u$ and $v$ which must equal to $ a_i + b
$ and $a_{i +1} +b $ with some $i\in J_A$ and $b\in B$,
then their difference is $ v - u = d\in D(A) $, so $ i= i(d)
$ where $ d = a_{i+1} - a_i $
which are thus uniquely determined and hence so
is $ b = u-a_i \ (=v-a_{i+1}) $, therefore there is, at most, one
such pair $ ( i,b) $. We
therefore deduce that
$$
P \leq \sum_{j=1}^k \binom{ |I_j|
} {2}
$$
where $ |I_j| $ denotes the number of elements of $ A + B $ in $ I_j $.
If $|A+B|=kL+r,\ 0\leq r<k$, then one can certainly arrange the end points of the arcs $ I_j $ so that $r$ of the $|I_j|$ are $ = L+1$, the other $k-r$ of the $|I_j|$ are $ =L$, in which case
$$
P \leq r \binom{L+1}{2} + (k-r) \binom{L}{2} = 
\frac{(kL+r)^2}{2k} - \frac{r^2}{2k} - \frac{kL}{2} \leq 
\frac{|A+B|^2}{2k}.
$$
Comparing this with (3) gives the upper bound
$$
|D(A)| \leq k + \frac{|A+B|^2}{2k|B|},
$$
which we minimize by selecting $ k $ to be the integer satisfying
$\frac{|A+B|}{\sqrt{2|B|}}+1 > k \geq \frac{|A+B|}{\sqrt{2|B|}} $. This implies the result.
\enddemo

\demo{Deduction of Corollary 2}
$S_\alpha(N)=\{\alpha n\,:\,1\leq n\leq N\}$ is an arithmetic progression, and $S_\alpha(N)+S_\alpha(N)=S_\alpha(2N)\setminus \{\alpha\}$, so we have $|S_\alpha(N)+S_\alpha(N)|=2N-1$. The result readily follows from
Corollary 1 (or taking $A=B=S_\alpha(N)$ in Theorem 1).
\enddemo

\demo{Proof of the second part of Theorem 1}
(Construction of $A$ with a lot of differences by the greedy algorithm.)
Let $ a_1=b_1,\ a_2=b_2 $. We select each $ a_j=b_{k_j}, j\geq 3 $ to
be the minimum $k> k_{j-1} $
such that $b_k - a_{j-1} $ is different
from the previous consecutive differences, $a_2 - a_1 , a_3 - a_2,...,a_{j-1}
- a_{j-2} $. Evidently we need to
avoid $ j-2 $ values so $ k_j \leq k_{j-1}+j-1 $. Hence
$$
k_j \leq (j-1) +
(j-2) +...+ 1 + 1 = \binom {j}{2} +1
$$
and there is room for choosing $k_j$ while this is $\leq |B| $.
Hence we can take
$\sqrt{2|B|} > j \geq \sqrt{2|B|}-1$. This gives the second
part of Theorem 1.
\enddemo

Another proof, though with a weaker constant, follows from recalling
that every set $B$ contains
a Sidon set $A$ of size $\gg \sqrt{|B|}$; and that
there are $|A|-1$ distinct differences between consecutive elements of
a Sidon set $A$.

\bigskip\head 4. The number of differences in the
original set \endhead

\demo{Proof of Theorem 2} Think about our sets in the  anticlockwise  ``circular order''. For any $b'$, $b'' \in B $ there exists $c\in C, b\in B $
such that $ c-b = b' -b'' $. Consider the anti-clockwise oriented arc from $b$ to $c$, this has the same arc length as the anti-clockwise oriented arc from $b''$ to $b'$. Locate the points $b=b_0,\dots, b_r=c$ of $B$ sitting on this arc and ordered in anti-clockwise direction. Each $ b_{i+1} - b_i$ is of the form $ c^* - b^* $, and  we can keep dividing those arcs into
smaller and smaller parts until they can get no smaller. But then the
parts must be of the
form $ c - b_c^- $ where there are no elements in--between
$ c $ and $ b_c^- $. Hence, $ b_c^-$ must
be the {\sl smaller neighbour} of $c$. This proves the first half of the theorem. The other half comes from the exact same analysis with choosing clockwise oriented arcs at the very beginning and always after. Finally, to see that $ \langle R^- \rangle_{\Bbb N_0} = \langle R^+ \rangle_{\Bbb N_0}$ note that both $R^-\subset B-B \subset  \langle R^+ \rangle_{\Bbb N_0} $ and $R^+\subset B-B \subset \langle R^- \rangle_{\Bbb N_0}$. 
\enddemo

\demo{Proof of Proposition 1} Let $S$ be a set of positive integers $ \leq N$, which contains no non--trivial 3--term arithmetic progressions. 
That is, if $2s=s_1+s_2$ in $S$, then $s_1=s_2=s$. There is such a set of $r_3(N)$ elements. Select $x=5N-2|S|$. Let
$$
B = S \cup \{ 2N< n \leq x -2 N \} \cup (x - S) ,
$$
so that $ |B| = x-4N + 2 |S|=N $ and $B + B \subseteq \{1,2,...,2x\}$, so that $|B+B|\leq 2x<10N$.
Now  for each $s\in S$ consider $ m = x-2s = (x-s) - s \in B-B $. Note that $m\geq x-2N$. 
If for some $b_1,b_2\in B$ we have $ m= b_1-b_2 $ then  
$ b_1 > x-2N$, else $ b_2= b_1-m \leq x-2N - (x-2N) $ which is impossible,  so $b_1=x-s_1,\ s_1\in S$. 
Similarly  $b_2 \leq 2N$ else $b_1=b_2+m > 2N+(x-2N)$ which is impossible,   and so $b_2=s_2,\ s_2\in S$. 
Therefore $2s=x-m=x-b_1+b_2=s_1+s_2$, which is only possible in $S$ when $s_1=s_2=s$. 
Hence if $C-B=B-B$ then  $x-S\subset C$, and so $|C|\geq |S|=r_3(N)$.  
\enddemo

\bigskip\head 5. Nearest neighbours \endhead

In this section we investigate the set $\{ N_a-a:\ a\in A\}$ of vectors from points of
$A$ to their nearest neighbours $N_a$ in $A$.

We begin with an important lemma from the theory of sphere packing:

\proclaim{Lemma 1} There exists an absolute constant $\kappa>0$ such that if
$z_1, z_2,\ldots, z_k\in \Bbb R^d$ with $|z_i|=1$ and $|z_i-z_j|\geq 1$ for all
$i\neq j$ then $k\leq \kappa (4/3)^d$.
\endproclaim

This is the celebrated ``{\sl kissing problem}'': how
many points can one place on the surface of a unit sphere such that they are all separated by
Euclidean distance at least 1? The exact answer is known only for $d = 1, 2, 3, 4, 8,$
and $24$. According to
[4, pages 23-24], the best bounds known in general are $1.15^{d}\ll k\leq \kappa \cdot 1.33^{d}$.

\proclaim{Lemma 2} There exists an absolute constant $\kappa>0$ such that if
$z_1, z_2,\ldots, z_k\in (\Bbb R/ \Bbb Z)^d$ with $\| z_i-z_j\| \geq \max \{\| z_i\| ,\| z_j\| \}$ for all
$i\neq j$ then $k\leq \kappa (4/3)^d$.
\endproclaim

\demo{Proof} For each $z_i\in (\Bbb R/\Bbb Z)^d$ we choose that   representative,  $\tilde {z_i} \in \Bbb R^d$ which is closest to the origin. Draw a ``radial half line'' from the origin across each $\tilde {z_i}$. Note that $\| z_i\| $ is precisely the Euclidean length of $\tilde {z_i}$. Moreover, $| \tilde {z_i} - \tilde {z_j}|  \geq \| z_i - z_j\| $, because this latter one is the minimal distance between any representatives of the two points. Thus we have
$$|\tilde {z_i} - \tilde {z_j}| \geq \max \{|\tilde {z_i}|,|\tilde {z_j}|\} \text{ for all } i\neq j.\tag 4$$

Now we have $k$ half lines, starting from the origin and each supporting exactly one point $\tilde {z_i}$. Two points on one half line violates (4). Even more, any two such radial lines must be at an angle of at least $\pi/3$, a smaller angle also violates (4).
If $r=\max_j |\tilde {z_j}|$ then we can move each point out to the circumference
of the sphere of radius $r$ along its ``radial half line''. This spreads the points out as far as possible since the angle between any two radial lines is at least $\pi/3$. A renormalized Lemma 2 proves the statement.
\enddemo

\proclaim{Corollary 3} There exists an absolute constant $\kappa>0$ such that
any $z\in (\Bbb R/\Bbb Z)^d$ belongs to $\leq \kappa (4/3)^d$ balls of the form
$B_a(\| N_a-a\| )$ (where $B_x(r)$ denotes the closed ball of radius $r$ centered at $x$.)
\endproclaim

\demo{Proof}  Suppose that $z\in B_{a_j}(\| N_{a_j}-a_j\| )$ for $j=1,2,\dots ,k$.
Then $\| a_i-a_j\| \geq \| N_{a_j}-a_j\|  \geq \| z-a_j\| $ whenever $i\ne j$. Writing
$a'_j=z-a_j$ for all $j$ we find that $\| a_i'-a_j'\| =\| a_i - a_j\| \geq \| a_j'\| $ for all $i\ne j$, and so
$k< \kappa (4/3)^d$ by Lemma 2.
\enddemo

\demo{Proof of Theorem 3}  Note that $w_i-w_j = w_{i-j} \in (\Bbb R/\Bbb Z)^d$, and $w_{-i} = -w_i$, hence
$D$ is a subset of  $\{\pm w_1, \pm w_2,\ldots, \pm w_{N-1}\}$, which we re-order so
that $\| w_{k_1}\| \leq \| w_{k_2}\|  \leq \ldots $. Recall that $j(i)$ minimizes
$\| w_i - w_{j(i)}\| $. If $j=i+k$ or $j=i-k$ then $\| w_i - w_j\| =\| w_k\| $, so to find  $j(i)$ we have to find the smallest integer $\ell$ such that either $i+k_\ell$ or $i-k_\ell$ falls into the interval $[1,N]$. 
Observe that if $k\leq\frac N2$ then 
$i+k$ is a legal choice for all $i\leq \frac N2$, and $i-k$ is  for all $i> \frac N2$. Let $\ell$ be the smallest positive integer such that $k_\ell\leq \frac N2$. That means $k_1,\dots,k_{\ell-1}>\frac N2$, and as $k_\ell$ is a legal choice for every $i$, we have restricted the choices to $D\subset
\{\pm w_{k_1}, \pm w_{k_2},\ldots, \pm w_{k_\ell}\}$.

Now if $1\leq u<v<\ell$ then $\frac N2 < k_u,k_v\leq N$ so that
$j=| k_u-k_v| < \frac N2$ but $k_\ell$ is the first element
of the reordered sequence which is $\leq \frac N2$, so $j\geq \ell$. Hence $\| w_{k_u}-w_{k_v}\| =\| w_j\| \geq \| w_{k_\ell}\| 
> \max\{ \| w_{k_u}\| , \| w_{k_v}\| \}$, and the result follows from Lemma 2.
\enddemo

\demo{Proof of Theorem 4}
Let $B\neq \emptyset$ be any fixed finite subset of $(\Bbb R/\Bbb Z)^d$.
For any $a\in A$ the ball $B_a(\| N_a - a\| )$ contains no elements of $A$ other than $a$
itself at the center, $N_a$ and possibly some others on the boundary.
If we translate this ball by an element $b\in B$ then
$B_{a+b}(\| N_a - a\| )$ may contain many elements of $A+B$, but
usually it is not the case. To see this, fix an element $b\in B$, and consider
$$
\align
\sum_{a\in A} | (A+B) \cap B_{a+b}(\| N_a-a\| ) | &= 
\sum_{c\in A+B} \#\{ a\in A: c\in B_{a+b}(\| N_a-a\| )\} \\
&= \sum_{c\in A+B} \#\{ a\in A: c-b\in B_a(\| N_a-a\| )\}  \leq
|A+B|\kappa(4/3)^d,\tag 5
\endalign
$$
because Corollary 3 says that $c-b$ cannot be in too many balls of the form $B_a(\| N_a-a\| )$.
On the other hand, let $\Upsilon_b$ be the set of $a\in A$ for which
$$
\left| (A+B) \cap B_{a+b}( \| N_a-a\| ) \right| > \frac{2\kappa}{\epsilon} (4/3)^d \frac{|A+B|}{|A|}.
$$
Comparing with (5) we have
$$
|\Upsilon_b| \frac{2\kappa}{\epsilon} (4/3)^d \frac{|A+B|}{|A|} <
\sum_{a\in \Upsilon_b} | (A+B) \cap B_{a+b}(\| N_a-a\| ) | \leq
|A+B|\kappa(4/3)^d,
$$
and hence $|\Upsilon_b|<\frac\epsilon2 |A|$.

Define 
$$A_c:=\{a\in A: c-a=b  \ \text{for some} \ b\in B, \ \text{and}\ a\not\in \Upsilon_b\},$$
 and
$N_c^{(l)}$ be the set of $l$ nearest neighbours to $c$ in $A+B$, where
$l=1+[\frac{2\kappa}{\epsilon} (4/3)^d  |A+B|/|A|]$.  If $a\in A_c$ and $c=a+b$ then $a\not\in \Upsilon_b$ so,
$B_c(\| N_a-a\| )$ contains $<l$ elements of $A+B$, though it does contain $N_a+b$, and therefore $N_a+b\in N_c^{(l)}$. Hence
$$
N_a-a = (N_a+b)-c\in N_c^{(l)}-c.
$$
This means that there are $\leq l$ choices for $N_a-a$,
as $a$ runs over all of the elements of $A_c$.

Let $R_1=\emptyset$.
Select $c_1  \in A+B$ so that $A_{c_1}\setminus R_1$ is of maximal size, and then let $R_2=A_{c_1}\cup R_1$. More generally, given $R_{k}$ we select  $c_k  \in A+B$ so that $A_{c_k}\setminus R_k$ is of maximal size, and then let $R_{k+1}=A_{c_k}\cup R_k =\bigcup_{j=1}^{k} A_{c_j}$, for each $k\geq 1$. Therefore
$$
\align
\sum_{c\in A+B} |A_c\setminus R_k| &=
\sum_{b\in B} \#\{ a\in A\setminus R_k: a\not\in \Upsilon_b\} \\
&\geq \sum_{b\in B} (|A|-|R_k|-|\Upsilon_b|) \\
&\geq |B|( (1-\frac{\epsilon}2)|A|-|R_k|).
\endalign
$$
Hence, by definition,  $A_{c_k}$ satisfies
$$
|A_{c_k}\setminus R_k|  \geq \frac{1}{|A+B|} \sum_{c\in A+B} |A_c\setminus R_k| \geq 
\frac{|A||B|}{|A+B|}
\left( 1-\frac{\epsilon}2 - \frac{|R_k|}{|A|} \right) .\tag 4
$$
We define $\theta_j\in [0,1]$ by $|R_j|=(1-\theta_j)|A|$, so that $(\theta_j)_{j\geq 1}$ is a decreasing sequence.
By (4) we have
$$
(\theta_k-\theta_{k+1})|A|=|R_{k+1}| -|R_k| = |A_{c_k}\setminus R_k| \geq 
\frac{|A||B|}{|A+B|} \big(\theta_k-\frac{\epsilon}2 \big),
$$
Therefore if $\theta_k\geq \epsilon$   then
$$
\theta_{k+1}\leq \left( 1-\frac{|B|}{2|A+B|}\right) \theta_k.
$$
We deduce that $\theta_n< \epsilon$ for some $n\ll  \frac{|A+B|}{|B|} \log(1/\epsilon)$, else
$\theta_k\geq \epsilon$ for all $k\leq n$ and so 
$$
\theta_{n}\leq \left( 1-\frac{|B|}{2|A+B|}\right)^{n-1}    < \epsilon.
$$
as  $\theta_1=1$.  

Finally we let  $A'=R_n$ so that $A'\subset A$ with $|A'|>(1-\epsilon)|A|$. Now $A'$ is the union of $n-1$ sets of the form $A_{c_j}$; for each $j$ we select $a_j\in A_{c_j}$ which maximizes $| (A+B) \cap B_{c_j}(\| N_{a_j}-a_j\| ) | $. Therefore
$$
|\{ N_a-a:\ a\in A'\}| \leq \sum_{j=1}^{n-1} | (A+B) \cap B_{c_j}(\| N_{a_j}-a_j\| ) | <
nl \ll (4/3)^d\cdot\frac{|A+B|^2}{|A||B|}
\cdot \frac1\epsilon\log\frac1\epsilon,
$$
and hence the result.
\enddemo

\example{Example} We construct an example 
$$
A=\{ 1, 4, 9, 16,\ldots ,m^2, m^2+1, m^2+2,\ldots , 2m^2-m\} \subset \Bbb R
$$
and convert the example to $\Bbb R/\Bbb Z$, by considering $\{ a/4m^2:\ a\in A\}\subset \Bbb R/\Bbb Z$.
Now $|A|=m^2, \  |A+A| < 4m^2$ and $|\{ N_a-a:\ a\in A\}|\geq m$. So if we
take $\epsilon =1/(2m)$ then, for any subset $A'\subset A$ with $|A'|>(1-\epsilon)|A|$. we have
$$
|\{ N_a-a:\ a\in A'\}|\geq |\{ N_a-a:\ a\in A\}|-|A\setminus A'|\geq m-\epsilon m^2\geq m/2
$$
whereas
$$
(4/3)^d  \cdot \min_B \frac{|A+B|^2}{|A||B|} \ \cdot \  \frac 1{\epsilon} \log \frac 1{\epsilon}
\leq (4/3) \frac{|A+A|^2}{|A|^2} \cdot 2m \log 2m \ll m \log 2m.
$$
This exhibits  that the upper bound in Theorem 4 cannot be much improved.
\endexample

\bigskip

\noindent{\sl Acknowledgements}: Thanks to Henry Cohn for some help with sphere packing,
and Boris Bukh and Seva Lev for their remarks.

\enddocument